\documentclass[12pt]{article}
\usepackage[utf8]{inputenc}
\usepackage{amsmath}
\usepackage{amssymb}
\usepackage{mathrsfs}
\usepackage{amsthm}
\usepackage{cite}
\usepackage{hyperref}
\usepackage{enumerate}
\usepackage[left=0.7in,right=0.7in,top=0.7in,bottom=0.7in]{geometry}
\usepackage{titlesec}
\titleformat*{\section}{\large\bfseries}
\newtheorem{theorem}{Theorem}[section]

\newtheorem{corollary}[theorem]{Corollary}
\newtheorem{definition}[theorem]{Definition}
\newtheorem{example}[theorem]{Example}
\newtheorem{proposition}[theorem]{Proposition}

\numberwithin{equation}{section}
\title{Product isometry of generalized weighted composition operator on general weighted Hardy space}
\author{\large Anuradha Gupta and Geeta Yadav$^*$}
\date{}
\begin{document}
\maketitle
\begin{abstract} 
We obtain necessary and sufficient conditions for the composition and weighted composition operator and product of composition operators to be isometry and unitary on $H_{E}(\xi).$ With the help of counter example we also prove that the product of two non isometric composition operator and weighted composition operator can be isometry on $H_{E}(\xi)$. We also completely characterize the boundedness of generalized weighted composition operators on $H_{E}(\xi).$ 


\end{abstract}    

\section{Introduction and preliminaries}

Let $\xi=(\xi_n)_{n\in \mathbb{Z_+}}$ be a sequence of positive real numbers satisfying 
\begin{equation} \label{condition on xi sequence}
lim_{n\rightarrow \infty} (\xi_n)^{1/n}= \infty
\end{equation}
where $\mathbb{Z_+}$ denotes the set of non negative integers.
By Hilbert space $H_E(\xi)$ (see \cite{Main Paper on H_E(beta)}), we denote the weighted Hardy space of entire functions induced by $\xi$ as
$$H_E(\xi)=\Big\{f(z)= \sum_{n=0}^{\infty} a_n z^n : ||f||^2=\sum_{n=0}^{\infty} |a_n|^2 {\xi_n}^2 <\infty \Big\}$$ 
and for $f(z)=\sum_{n=0}^{\infty} b_n z^n\,\,\text{and}\,\,g(z)=\sum_{n=0}^{\infty} c_n z^n \,\text{in}\, H_E(\xi)$ the inner product is defined as
\begin{equation} \label{inner product define}
\langle f,g \rangle=\sum_{n=0}^{\infty} b_n \overline{c_n} {\xi_n}^2. 
\end{equation}
Elements of Hilbert space $H_E(\xi)$ are entire functions if and only if \eqref{condition on xi sequence} holds \cite{Hilbert space of entire fun}.\\
For the sequence $\xi=(\sqrt{n!})_{n\in \mathbb{Z_+}},$ Hilbert space $H_E(\xi)$  is the well known Fock space $\mathcal{F}^2.$ The sequence $e_n(z)=\frac{1}{\xi_n}z^n$ for $n\in \mathbb{Z_+}$ forms an orthonormal basis for $H_E(\xi).$ \\

The space $H_E(\xi)$ is also reproducing kernel Hilbert space. The reproducing kernel, $K_{p}\in H_E(\xi),$ for point evaluation at $p \in \mathbb {C}$ satisfying $f(p)=\langle f, K_p \rangle$ for each $f\in H_E(\xi),$ is
$$K_p(z)= \sum_{n=0}^{\infty} \frac{1}{{\xi_n}^2} \overline{p}^n z^n. $$

Let $W$ be a subset of $\mathbb{C}$ and $H(W)$ denote the Hilbert space of analytic functions on $W.$ If $\Phi$ and $\Upsilon$ are holomorphic functions on $W$ such that $\Phi$ maps $W$ into $W$ then, for $f\in H(W)$ 
$$C_{\Upsilon, \Phi}f=\Upsilon\cdot (f \circ \Phi) $$
denote the weighted composition operator and if $\Upsilon\equiv 1,$ then it is the composition operator $C_{\Phi}f=f \circ \Phi.$ Researchers, find it interesting to analyze the relation between the operator- theoretic behaviours of $C_{\Upsilon, \Phi}$ and function-theoretic properties of $\Upsilon$ and $\Phi.$ These relations have been studied on various function spaces on the unit disk like Hardy space, Bergman space, derivative Hardy space and on Hilbert spaces of entire functions like Fock space. In this paper we discuss the results on $H_E(\xi),$ the weighted Hardy space of entire functions on $\mathbb{C}.$ \\   

Let $\Phi,\Upsilon:\mathbb{C} \longrightarrow \mathbb{C} $ be entire functions and $p$ be a nonnegative integer. The formal generalized weighted composition expression on $ H_E(\xi)$ is defined as
$$T(p,\Phi,\Upsilon)f=\Upsilon \cdot f^p \circ \Phi$$
where $f^p$ is the $p$th derivative of function $f$ in $ H_E(\xi)$. The generalized weighted composition operator $\mathcal{D}_{\Phi,\Upsilon}^p$ (see \cite{Zhu Generalized weighted 1},\cite{Zhu Generalized weighted 2}) on $H_E(\xi)$ corresponding to the  expression $T(p,\Phi,\Upsilon)$  is defined by  
  
$$ \mathcal{D}_{\Phi,\Upsilon}^p f= T(p,\Phi,\Upsilon)f \, \text{ for}\,f\in Dom( \mathcal{D}_{\Phi,\Upsilon}^p)$$ 
where $Dom(\mathcal{D}_{\Phi,\Upsilon}^p)=\{f \in H_E(\xi): T(p,\Phi,\Upsilon)f\in H_E(\xi)\}.$ For $p=0,$ $\Upsilon \equiv 1,$ $ \mathcal{D}_{\Phi,\Upsilon}^p$ is the composition operator $C_{\Phi}$ induced by $\Phi$; for $p=0,$  it is the weighted composition operator $C_{\Upsilon,\Phi}=\Upsilon C_{\Phi}$ and for  $\Upsilon \equiv 1, \Phi(z)=z$ it is $\mathcal{D}^p,$ the $p$-th differentiation operator.\\

Throughout the paper we assume that $\Phi$ and $\Upsilon$ are entire functions on $\mathbb{C}.$ A bounded linear operator $S$ on a Hilbert space $H$ is said to be isometry if $||Sx||=||x||$ for all $x\in H$ or equivalently $\langle Sx, Sy \rangle=\langle x, y \rangle$ for all $x,y \in H$ or equivalently $S^*S=I_{H}.$ In section 2, the necessary and sufficient condition for composition operator and product of composition operators to be isometry or unitary on $H_E(\xi)$ are obtained. In section 3, we study the results of weighted composition operators on $H_E(\xi)$. In section 4, we discuss  m-isometry and invertible composition operators on $H_E(\xi)$. In the last section, we give the characterization for boundedness of generalized weighted composition operators on $H_E(\xi)$.

\section{Isometry composition operator} 

We begin this section with the definition of the adjoint of a bounded linear operator on $H_E(\xi).$ 
\begin{definition}
The adjoint of a bounded linear operator $S$  on $H_E(\xi),$ denoted by $S^*$, is the operator which satisfies
$$\langle Sg, f \rangle=\langle g,S^*f \rangle $$
for $f,g \in H_E(\xi).$ 
\end{definition}

For function $\Phi$ defined by  $\Phi(z)=\mu z$ for $z\in \mathbb{C}$ and $|\mu|\leq 1,$ the composition operator induced by $\Phi$ is bounded on $H_E(\xi).$ 
\begin{theorem} \label{Com Op isometry iff mod(a)=1}
Let $|\mu|\leq 1$ and  $\Phi(z)=\mu z$ for $z\in \mathbb{C}.$ Then, $C_{\Phi}$ is an isometry on $H_E(\xi)$ if and only if $|\mu|=1.$ 
\end{theorem} 
\textbf{Proof.} Let $\Phi(z)=\mu z.$ Then,
\begin{equation} \label{adjoint of CPhi for Phi mz}
C_{\Phi}^*=C_{\Psi} \,\,\text{ where}\,\, \Psi(z)=\overline{\mu} z \,\, \text{for} \,\, z\in \mathbb{C}.
\end{equation} 
For $f\in H_E(\xi)$ and $z\in \mathbb{C}$ 
\begin{equation}\label{equation f(lambda lambda bar z)}
(C_{\Phi}^* C_{\Phi}f)(z)=(C_{\Psi}C_{\Phi}f)(z)=((f\circ \Phi)\circ \Psi)(z) =f(|\mu|^2 z). 
\end{equation}
Let $C_{\Phi}$ be an isometry. Then, for $f\in H_E(\xi)$ and $z\in \mathbb{C}$
\begin{equation} \label{equ 2 Cphi isometry}
(C_{\Phi}^* C_{\Phi}f)(z)=f(z).
\end{equation}
In particular, for function $f(z)=z$ in $H_E(\xi),$ equations (\ref{equation f(lambda lambda bar z)}) and (\ref{equ 2 Cphi isometry}) give $|\mu|^2 z=z.$ Then, for $z=1,$ $|\mu|=1.$ \\
Conversely, let $|\mu|=1.$ Then, $(C_{\Phi}^* C_{\Phi}f)(z)=f(z).$ This implies $C_{\Phi}^* C_{\Phi}=I_{H_E(\xi)}$ and hence, $C_{\Phi}$ is an isometry. \\

In the following result we characterize the unitary composition operator induced by $\Phi$ defined in Theorem \ref{Com Op isometry iff mod(a)=1}.
\begin{theorem} \label{Com Op unitary iff mod(a)=1}
Let $|\mu|\leq 1$ and  $\Phi(z)=\mu z$ for $z\in \mathbb{C}.$ Then, $C_{\Phi}$ is unitary on $H_E(\xi)$ if and only if $|\mu|=1.$ 
\end{theorem} 
\textbf{Proof.} Let $C_{\Phi}$ be unitary. Then, $C_{\Phi}$ is isometry and by Theorem \ref{Com Op isometry iff mod(a)=1}, it follows that $|\mu|=1.$\\
Conversely, let $|\mu|=1$ then by Theorem \ref{Com Op isometry iff mod(a)=1}, $C_{\Phi}^*C_{\Phi}=I_{H_E(\xi)}.$ To prove $C_{\Phi}$ is unitary, it is sufficient to prove that $C_{\Phi}C_{\Phi}^*=I_{H_E(\xi)}.$ Now using equation \eqref{adjoint of CPhi for Phi mz},  for $f\in H_E(\xi)$ and $z\in \mathbb{C}$
\begin{equation} \label{equation f(lambda bar lambda  z)}
(C_{\Phi}C_{\Phi}^*f)(z)=(C_{\Phi}C_{\Psi}f)(z)=((f\circ \Psi)\circ \Phi)(z)=f(|\mu|^2 z)=f(z)
\end{equation}
Thus $C_{\Phi}C_{\Phi}^*=I_{H_E(\xi)}.$\\

Combination of Theorem \ref{Com Op isometry iff mod(a)=1} and Theorem \ref{Com Op unitary iff mod(a)=1} lead to the following result:
\begin{corollary}
Let $|\mu|\leq 1$ and  $\Phi(z)=\mu z$ for $z\in \mathbb{C}.$ Then, $C_{\Phi}$ is isometry on $H_E(\xi)$ if and only if it is unitary.
\end{corollary}

\begin{theorem} \label{product lambda1z and lambda2z comp op}
Let $\Phi_1(z)=\mu_1 z$ and $\Phi_2(z)=\mu_2 z$ for $|\mu_1|\leq 1$ and $|\mu_2|\leq 1$ respectively. Then, the product $C_{\Phi_{1}}C_{\Phi_{2}}$ is isometry on $H_E(\xi)$ if and only if $|\mu_1|= 1$ and $|\mu_2|=1.$ Moreover, the product $C_{\Phi_{1}}C_{\Phi_{2}}$ is unitary on $H_E(\xi)$ if and only if $|\mu_1|= 1$ and $|\mu_2|=1.$
\end{theorem}
\textbf{Proof.} For $f\in H_E(\xi)$ and $z\in \mathbb{C}$ 
$$(C_{\Phi_1}C_{\Phi_2}f)(z)=((f \circ \Phi_2)\circ\Phi_1)(z)=f((\mu_2 \mu_1)z).$$
Thus, $C_{\Phi_1}C_{\Phi_2}=C_{\Phi}$ where $\Phi(z)=(\Phi_{2} \circ \Phi_1)(z)=(\mu_2 \mu_1)z$ for $z\in \mathbb{C}$ and 
$|\mu_2 \mu_1|=|\mu_2| |\mu_1|\leq 1.$\\
Let $C_{\Phi_1}C_{\Phi_2}$ be isometry. Then, by Theorem \ref{Com Op isometry iff mod(a)=1}
 \begin{equation} \label{mod lambda 1lambda 2 equal 1}
 |\mu_2| |\mu_1|=|\mu_2 \mu_1|=1
 \end{equation} 
It follows clearly that $|\mu_1| \neq 0$, $|\mu_2| \neq 0.$ Now, since $|\mu_1|\leq 1$ and $|\mu_2|\leq 1,$ therefore, equation \eqref{mod lambda 1lambda 2 equal 1} gives
$$1 \leq \frac{1}{|\mu_2|}=|\mu_1|\leq 1.$$
Thus, $|\mu_1|=1$ and $|\mu_2|=1.$\\
Conversely, let $|\mu_1|=1$ and $|\mu_2|=1.$ Then, $|\mu_2 \mu_1|=1$ and by Theorem \ref{Com Op unitary iff mod(a)=1}, $C_{\Phi_1}C_{\Phi_2}=C_{\Phi}$ is unitary and hence, an isometry.\\

If  $\mu \in \mathbb {C}$ such that $|\mu| \leq 1,$ then $|\overline{\mu} |\leq 1.$ This along with equation \eqref{adjoint of CPhi for Phi mz} and Theorem \ref{product lambda1z and lambda2z comp op} leads to the following result:
\begin{corollary}
Let $\Phi_1(z)=\mu_1 z$ and $\Phi_2(z)=\mu_2 z$ for $|\mu_1|\leq 1$ and $|\mu_2|\leq 1.$ Then $C_{\Phi_{1}}^*C_{\Phi_{2}}$ $(\text{or}\,\, C_{\Phi_{1}}C_{\Phi_{2}}^* \text{or}\,\, C_{\Phi_{1}}^*C_{\Phi_{2}}^*)$ is isometry if and only if $|\mu_1|= 1$ and $|\mu_2|=1.$ Moreover, $C_{\Phi_{1}}^*C_{\Phi_{2}}$ $($or $ C_{\Phi_{1}}C_{\Phi_{2}}^*$ or $C_{\Phi_{1}}^*C_{\Phi_{2}}^*)$ is unitary if and only if $|\mu_1|= 1$ and $|\mu_2|=1.$ 
\end{corollary}
In the following result $\Phi$ (for $N=1$) is choosen according to the result on necessary condition given by  Doan et al. (\cite{Com Op on several variable}, Theorem 4.6) for composition operator $C_{\Phi}$ to be bounded on $H_E(\xi).$ In \cite{Main Paper on H_E(beta)} Tan and Khoi discussed the  boundedness of composition operators and obtained many important results on $H_E(\xi).$ The necessary and sufficient condition for a bounded composition operator on $H_E(\xi)$ to be isometry or unitary is as follows:
\begin{theorem} \label{isometry condition for az+b}
Let $\Phi(z)=\mu z+d$ for $|\mu|\leq 1$ and $d\in \mathbb{C}$ such that $C_{\Phi}$ is bounded on $H_E(\xi).$  Then, $C_{\Phi}$ is an isometry on $H_E(\xi)$ if and only if $|\mu|= 1$ and $d=0.$ Moreover, $C_{\Phi}$ is unitary on $H_E(\xi)$ if and only if $|\mu|= 1$ and $d=0.$
\end{theorem}
\textbf{Proof.} Let $C_{\Phi}$ be an isometry. Then, for $g,h\in H_E(\xi)$
\begin{equation*} 
 \langle C_{\Phi}h, C_{\Phi}g \rangle=\langle h, g\rangle 
\end{equation*}
that is,
\begin{equation} \label{isometry az+b condition 1 equ}
\langle h\circ \Phi, g\circ \Phi \rangle=\langle h, g\rangle.
\end{equation} 
By our assumption $C_{\Phi}$ is isometry implies, $d=0$ because if $d\neq 0$ then for functions $h(z)=1$ and $g(z)=z\,\, \text{for all}\,\, z\in \mathbb{C}$ in $H_E(\xi)$
\begin{equation}  \label{isometry az+b condition 1* equ}
(C_{\Phi}h)(z)=1 \,\,\text{and} \,\, (C_{\Phi}g)(z)=\Phi(z)=\mu z+d \,\,\,\, \text{for all}\, \,z\in \mathbb{C}.
\end{equation}
Since $\xi_0$ is positive, therefore, by combining equations (\ref{inner product define}) and (\ref{isometry az+b condition 1* equ}) we get
$$\langle h\circ \Phi, g\circ \Phi \rangle=\overline{d}{\xi_0}^2 \neq 0=\langle h, g\rangle$$ 
which contradicts equation (\ref{isometry az+b condition 1 equ}). Hence $\Phi(z)=\mu z$ for all $z\in \mathbb{C}.$\\
Now, if  $h(z)=z$ and $g(z)=z$ for all $z\in \mathbb{C}$ in $H_E(\xi)$ then 
\begin{equation}  \label{isometry az+b condition 1** equ}
(C_{\Phi}h)(z)=\Phi(z)=\mu z \,\,\text{and} \,\, (C_{\Phi}g)(z)=\mu z \,\,\,\, \text{for all} \,z\in \mathbb{C}.
\end{equation}
Combining equations (\ref{inner product define}), (\ref{isometry az+b condition 1 equ}) and (\ref{isometry az+b condition 1** equ}) we obtain
$$|\mu|^2 {\xi_1}^2={\xi_1}^2$$
or
\begin{equation} \label{isometry az+b condition 2 equ}
(|\mu|^2-1){\xi_1}^2=0.
\end{equation}
Since ${\xi_1}$ is positive, therefore, from equation (\ref{isometry az+b condition 2 equ}) we get $|\mu|=1.$\\
Conversely, let $|\mu|=1$ and $d=0$ then, by Theorem \ref{Com Op unitary iff mod(a)=1} it follows that $C_{\Phi}$ is unitary and hence isometry.
\begin{theorem} \label{product comp op az+b is isometry }
Let $\Phi_1(z)=\mu_1 z+d_1\,\,,|\mu_1|\leq 1$ and $\Phi_2(z)=\mu_2 z+d_2 \,\,,|\mu_2|\leq 1$ such that $C_{\Phi_1}$ and    $C_{\Phi_2}$ are bounded on $H_E(\xi).$ Then, the product $ C_{\Phi_2}C_{\Phi_1}$ is isometry on $H_E(\xi)$ if and only if $|\mu_1|=1$, $|\mu_2|=1$ and $\mu_1d_2+d_1=0.$ Moreover, the product $ C_{\Phi_2}C_{\Phi_1}$ is unitary on $H_E(\xi)$ if and only if $|\mu_1|=1$, $|\mu_2|=1$ and $\mu_1d_2+d_1=0.$
\end{theorem} 
\textbf{Proof.}  Clearly $C_{\Phi_2}C_{\Phi_1}=C_{\Phi_1 \circ \Phi_2},\,\,|\mu_1 \mu_2|=|\mu_1| | \mu_2|\leq 1$  and for $z\in \mathbb{C}$
$$(\Phi_1\circ \Phi_2)(z)=\Phi_1(\Phi_2(z))=(\mu_1 \mu_2)z+(\mu_1d_2+d_1).$$  
Let $C_{\Phi_2}C_{\Phi_1}$ be isometry. Then, by Theorem \ref{isometry condition for az+b}, $|\mu_1 \mu_2|=1$ and $\mu_1d_2+d_1=0.$  Further, $|\mu_1|=1$, $|\mu_2|=1$ since $|\mu_1| |\mu_2|=|\mu_1 \mu_2|=1$ and $|\mu_1|\leq 1,$  $|\mu_2|\leq 1.$ \\
Conversely, let $|\mu_1|=1$, $|\mu_2|=1$ and $\mu_1d_2+d_1=0.$  Since $|\mu_1 \mu_2|=1$ by Theorem \ref{isometry condition for az+b}, $C_{\Phi_2}C_{\Phi_1}$ is unitary and hence isometry.\\
The following example shows that in Theorem \ref{product comp op az+b is isometry }, $C_{\Phi_2}C_{\Phi_1}$ can be isometry even when both $C_{\Phi_2}$ and $ C_{\Phi_1}$ are not isometry.
\begin{example}
Let $\Phi_1(z)=z+d$ and $\Phi_2(z)=z-d$ for all $z\in \mathbb{C}$ where $0\neq d\in \mathbb {C}$ such that $C_{\Phi_1}$ and    $C_{\Phi_2}$ are bounded on $H_E(\xi).$ Then, by Theorem \ref{isometry condition for az+b}, the composition operator $C_{\Phi_1}$ and $C_{\Phi_2}$ are not isometry on $H_E(\xi)$. Clearly,
$$(\Phi_1\circ \Phi_2)(z)=\Phi_1(z-d)=(z-d)+d=z\,\,\,\, \text{for all} \,\,z\in \mathbb{C}.$$
By Theorem \ref{product comp op az+b is isometry }, it follows that $C_{\Phi_2}C_{\Phi_1}$ is unitary and hence isometry on $H_E(\xi).$ 
\end{example}
The sum of two isometric composition operator $C_{\Phi_{1}}$ and $C_{\Phi_{2}}$ on $H_E(\xi)$ need not be an isometry. This can be verified from the following example.
\begin{example}
 Let $\Phi_1(z)=z$ and $\Phi_2(z)=-z$ for all $z\in \mathbb{C}.$ Then, by Theorem \ref{Com Op isometry iff mod(a)=1}, $C_{\Phi_{1}}$ and $C_{\Phi_{2}}$ are isometry on $H_E(\xi).$ \\
Now, for all $f\in \mathcal{F}^2$ and $z\in \mathbb{C}$ we have 
 \begin{align*}
((C_{\Phi_{1}}+C_{\Phi_{2}})f)(z)&=(f \circ \Phi_{1} +f \circ \Phi_{2})(z)\\
                                 &=f(\Phi_{1}(z))+f(\Phi_2(z))\\                                                                                 
                                 &=f(z)+f(-z).                                          
\end{align*}
In particular, for $g \in H_E(\xi)$ defined by $g(z)=z$ for all $z\in \mathbb{C}$ above equation reduces to 
$$((C_{\Phi_{1}}+C_{\Phi_{2}})g)(z)=g(z)+g(-z)=z+(-z)=0\,\, \text{for all}\,\, z\in \mathbb{C} $$
that is, $(C_{\Phi_{1}}+C_{\Phi_{2}})g$ is zero function.\\
Thus, $||(C_{\Phi_{1}}+C_{\Phi_{2}})g||=0\neq ||g||.$ Hence,  $C_{\Phi_{1}}+C_{\Phi_{2}}$ is not an isometry.
\end{example} 

\section{Isometric weighted composition operator} 
In this section, first we will discuss that for $\mu(\neq 0)$ and $\nu \in \mathbb{C}$ whether the sequence $\xi=(\xi_n)_{n\in \mathbb{Z_+}}$ satisfies the following equation or not.
\begin{equation} \label{eq for weighted comp isometry written before theorem}
{\xi_0}^2=|\nu|^{2n} |\mu|^{n(n-1)} {\xi_n}^2 \,\,\, \text{for}\,\, n\in \mathbb{N}.
\end{equation}
For $\xi=(\sqrt{n!})_{n\in \mathbb{Z_+}}$ and $n=1,2,3$  equation (\ref{eq for weighted comp isometry written before theorem}) reduces to 
\begin{equation} \label{eq 2 for weighted comp isometry written before theorem}
1=|\nu|^2 (1!), \,\, 1=|\nu|^4|\mu|^2 (2!) \,\, \text{and} \,\, 1=|\nu|^6|\mu|^6 (3!),  
\end{equation}
respectively. Infact, $|\nu|=1$ and substituting this value in the remaining two equalities of equation (\ref{eq 2 for weighted comp isometry written before theorem})  we obtain $|\mu|^2=\frac{1}{2}$ and $|\mu|^6=\frac{1}{6}$ which cannot be true for $|\mu|.$ Hence, sequence $\xi=(\xi_n)_{n\in \mathbb{Z_+}}$ can not satisfy the equation (\ref{eq for weighted comp isometry written before theorem}).\\

In 2014, Le \cite{Normal and isometric} obtained that an isometric weighted composition operator is unitary operator on $\mathcal{F}^2.$ The following result gives characterization for a weighted composition operator to be isometry on $H_{E}(\xi).$

\begin{theorem} \label{Weigh Com Op isometry on Hilbert sp Entire func iff}
Let $\Upsilon(z)=\nu z+c$ and  $\Phi(z)=\mu z$ for $z\in \mathbb{C}$ where $\nu, c,\mu \in \mathbb{C}$ with $|\mu|\leq 1$ such that $C_{\Upsilon,\Phi}$ is bounded on $H_E(\xi).$ Then $C_{\Upsilon,\Phi}$ is an isometry on $H_{E}(\xi)$ if and only if $|\mu|=1, \nu=0$ and $|c|=1.$ Moreover, $C_{\Upsilon,\Phi}$ is unitary on $H_{E}(\xi)$ if and only if $|\mu|=1, \nu=0$ and $|c|=1.$
\end{theorem} 
\textbf{Proof.} Let $C_{\Upsilon,\Phi}$ be an isometry. Then, $\mu\neq 0$ and for $g,h\in H_E(\xi)$
\begin{equation*} 
 \langle C_{\Upsilon,\Phi}h, C_{\Upsilon,\Phi}g \rangle=\langle h, g\rangle 
\end{equation*}
that is,
\begin{equation} \label{weigh comp isometry az+b condition 1 equ}
\langle \Upsilon\cdot (h \circ \Phi), \Upsilon\cdot (g \circ \Phi) \rangle=\langle h, g\rangle.
\end{equation} 
Now we consider the following cases for $h,g \in H_E(\xi)$ :\\
Case I. For $n\in \mathbb{Z_{+}},$ if $h(z)=z^n$ and $g(z)=z^n\,\, \text{for all}\,\, z\in \mathbb{C}$ then, 
\begin{equation}  \label{weigh comp isometry az+b powen n case equ}
(\Upsilon\cdot (h \circ \Phi))(z)=\Upsilon(z)h(\mu z)=\nu\mu^n z^{n+1}+c\mu^n z^n \,\,\text{and} \,\, (\Upsilon\cdot (g \circ \Phi))(z)=\nu\mu^n z^{n+1}+c\mu^n z^n.
\end{equation}
Now, combining equations (\ref{weigh comp isometry az+b condition 1 equ}) and (\ref{weigh comp isometry az+b powen n case equ}) we obtain
$$|\nu|^2 |\mu|^{2n}{\xi_{n+1}}^2+|c|^2|\mu|^{2n} {\xi_n}^2={\xi_n}^2$$
or
\begin{equation} \label{weigh comp isometry az+b condition n* equ}
|\nu|^2 |\mu|^{2n}{\xi_{n+1}}^2+(|c|^2|\mu|^{2n}-1) {\xi_n}^2=0.
\end{equation}
From equation (\ref{weigh comp isometry az+b condition n* equ}) it follows that if $h(z)=z^{0}=1$ and $g(z)=z^{0}=1\,\, \text{for all}\,\, z\in \mathbb{C}$ then 
\begin{equation} \label{weigh comp isometry az+b condition 2* equ}
|\nu|^2{\xi_1}^2+(|c|^2-1) {\xi_0}^2=0.
\end{equation}
Moreover, if $h(z)=z$ and $g(z)=z\,\, \text{for all}\,\, z\in \mathbb{C}$ then, 
\begin{equation} \label{weigh comp isometry az+b condition 2** equ}
|\nu|^2|\mu|^2 {\xi_2}^2+(|c|^2|\mu|^2-1) {\xi_1}^2=0.
\end{equation}
Case II. For $n\in \mathbb{Z_{+}},$ if $h(z)=z^n$ and $g(z)=z^{n+1}\,\, \text{for all}\,\, z\in \mathbb{C}$ then, 
\begin{equation}  \label{weigh comp isometry az+b powen odd even case equ}
(\Upsilon\cdot (h \circ \Phi))(z)=\Upsilon(z)h(\mu z)=\nu\mu^n z^{n+1}+c\mu^n z^n \,\,\text{and} \,\, (\Upsilon\cdot (g \circ \Phi))(z)=\nu\mu^{n+1} z^{n+2}+c\mu^{n+1}z^{n+1}.
\end{equation}
Now, combining equations (\ref{weigh comp isometry az+b condition 1 equ}) and (\ref{weigh comp isometry az+b powen odd even case equ}) we obtain
\begin{equation} \label{weigh comp isometry az+b condition odd even inner product equ}
\nu \bar{c} \mu^{n} \overline{\mu^{n+1}} {\xi_{n+1}}^2=0.
\end{equation}
Since $\xi_{n+1}\neq 0$ and $\mu \neq 0,$ therefore, from equation (\ref{weigh comp isometry az+b condition odd even inner product equ}) it follows that either $\nu=0$ or $c=0.$\\
Case(a) If $\nu=0$ then equations (\ref{weigh comp isometry az+b condition 2* equ}) and (\ref{weigh comp isometry az+b condition 2** equ}) give 
\begin{equation} \label{weigh Comp isometry case a is zero 1}
(|c|^2-1) {\xi_0}^2=0
\end{equation} 
and
\begin{equation} \label{weigh Comp isometry case a is zero 2}
(|c|^2|\mu|^2-1) {\xi_1}^2=0
\end{equation}
respectively. Since $\xi_0\neq 0$ from equation (\ref{weigh Comp isometry case a is zero 1}) it follows that $|c|^2=1$ and further substituting this value in equation (\ref{weigh Comp isometry case a is zero 2}) we get $|\mu|=1$ because $\xi_1\neq 0.$ \\
Case(b) If $c=0$ then equations (\ref{weigh comp isometry az+b condition 2* equ}) and (\ref{weigh comp isometry az+b condition 2** equ}) give
\begin{equation} \label{weigh Comp isometry case b is zero 1}
{\xi_0}^2=|\nu|^2{\xi_1}^2
\end{equation} 
and
\begin{equation} \label{weigh Comp isometry case b is zero 2}
 {\xi_1}^2=|\nu|^2|\mu|^2 {\xi_2}^2
\end{equation}
respectively. Combining equations (\ref{weigh Comp isometry case b is zero 1}) and (\ref{weigh Comp isometry case b is zero 2}) we get
$${\xi_0}^2=|\nu|^2{\xi_1}^2=|\nu|^4|\mu|^2 {\xi_2}^2.$$
With the help of equation (\ref{weigh comp isometry az+b condition n* equ}) it can be proved using Principle of Mathematical Induction that 
\begin{equation} \label{equation involving relation on xi0 and xin elements} 
{\xi_0}^2=|\nu|^{2n} |\mu|^{n(n-1)} {\xi_n}^2 \,\,\, \text{for}\,\, n\in \mathbb{N}.
\end{equation} 
But this can not be true for sequence $\xi=(\xi_n)_{n\in \mathbb{Z_+}}$. Hence, Case(b) is not possible.\\ 
Conversely, let $|\mu|=1, \nu=0$ and $|c|=1.$ Then, 
$$(C_{\Upsilon,\Phi}h)(z)=(\Upsilon\cdot (h \circ \Phi))(z)=\Upsilon(z)(C_{\Phi}h)(z)=c \,(C_{\Phi}h)(z)=((c\,C_{\Phi})h)(z).$$
Thus, $C_{\Upsilon,\Phi}=c\,C_{\Phi}.$ By our assumption $|\mu|=1,$ Theorem \ref{Com Op unitary iff mod(a)=1} implies that  $C_{\Phi}$ is unitary and since $|c|=1,$ therefore,
$$C_{\Upsilon,\Phi}^* C_{\Upsilon,\Phi}=(c\,C_{\Phi})^* (c\,C_{\Phi})=\overline{c}c C_{\Phi}^*C_{\Phi}=|c|^2I_{H_E(\xi)}=I_{H_E(\xi)}$$
and $$C_{\Upsilon,\Phi} C_{\Upsilon,\Phi}^*=(c\,C_{\Phi}) (c\,C_{\Phi})^*=c \overline{c}C_{\Phi}C_{\Phi}^*=|c|^2I_{H_E(\xi)}=I_{H_E(\xi)}.$$
Thus, $C_{\Upsilon,\Phi}$ is unitary and hence an isometry. 
\begin{theorem} \label{product of isometry weight comp op is isometry} 
 Let $\Upsilon_{1}(z)=\nu_{1}z+c_{1},\Upsilon_{2}(z)=\nu_{2}z+c_{2}, $  $\Phi_{1}(z)=\mu_{1} z$ and $ \Phi_{2}(z)=\mu_{2} z$ for $z\in \mathbb{C}$ where $\nu_{1},\nu_{2},c_{1},c_{2},\mu_{1},\mu_{2} \in \mathbb{C}$ with $|\mu_{1}|\leq 1$ and $|\mu_{2}|\leq 1,$ respectively, such that $C_{\Upsilon_{1},\Phi_{1}}$ and $C_{\Upsilon_{2},\Phi_{2}}$ are bounded on $H_E(\xi).$ Then, the product $C_{\Upsilon_{1},\Phi_{1}}C_{\Upsilon_{2},\Phi_{2}}$ of isometric weighted composition operators $C_{\Upsilon_{1},\Phi_{1}}$ and $C_{\Upsilon_{2},\Phi_{2}}$ is an isometry on $H_E(\xi).$ Moreover, the product $C_{\Upsilon_{1},\Phi_{1}}C_{\Upsilon_{2},\Phi_{2}}$ of isometric weighted composition operators $C_{\Upsilon_{1},\Phi_{1}}$ and $C_{\Upsilon_{1},\Phi_{1}}$ is unitary on $H_E(\xi).$
\end{theorem}
\textbf{Proof.} Let $C_{\Upsilon_{1},\Phi_{1}}$ and $C_{\Upsilon_{2},\Phi_{2}}$ be isometries on $H_E(\xi).$ Then, by Theorem \ref{Weigh Com Op isometry on Hilbert sp Entire func iff} $C_{\Upsilon_{1},\Phi_{1}}$ and $C_{\Upsilon_{2},\Phi_{2}}$ are unitaries. Since the product of two unitary is unitary, therefore, $C_{\Upsilon_{1},\Phi_{1}}C_{\Upsilon_{2},\Phi_{2}}$ is unitary and hence an isometry on $H_E(\xi).$\\

The following result follows from Theorem \ref{product of isometry weight comp op is isometry} on the fact that adjoint of a unitary operator is unitary.
\begin{corollary}
 Let $\Upsilon_{1}(z)=\nu_{1}z+c_{1},\Upsilon_{2}(z)=\nu_{2}z+c_{2}, $  $\Phi_{1}(z)=\mu_{1} z$ and $ \Phi_{2}(z)=\mu_{2} z$ for $z\in \mathbb{C}$ where $\nu_{1},\nu_{2},c_{1},c_{2},\mu_{1},\mu_{2} \in \mathbb{C}$ with $|\mu_{1}|\leq 1$ and $|\mu_{2}|\leq 1,$ respectively, such that $C_{\Upsilon_{1},\Phi_{1}}$ and $C_{\Upsilon_{2},\Phi_{2}}$ are isometries on $H_E(\xi)$. Then, $C_{\Upsilon_{1},\Phi_{1}}^*C_{\Upsilon_{2},\Phi_{2}}^*, C_{\Upsilon_{1},\Phi_{1}}^*C_{\Upsilon_{2},\Phi_{2}}$ and $C_{\Upsilon_{1},\Phi_{1}}C_{\Upsilon_{2},\Phi_{2}}^*$ are unitaries on $H_E(\xi)$. 
\end{corollary}

The converse of the Theorem \ref{product of isometry weight comp op is isometry} is not true. This can be verified from the following example:

\begin{example} \label{example product of non isometry is isometry}
 Let $\Phi_{1}(z)=iz,$ $\Phi_{2}(z)=-iz,$ $\Upsilon_{1}(z)=\frac{1}{2}$ and $ \Upsilon_{2}(z)=2$ for all $z\in \mathbb{C}.$ Then, by Theorem \ref{Weigh Com Op isometry on Hilbert sp Entire func iff}, $C_{\Upsilon_{1},\Phi_{1}}$ and $C_{\Upsilon_{2},\Phi_{2}}$ are not isometries on $H_{E}(\xi)$ because $\Upsilon_{1}$ and $\Upsilon_{2}$ are not constant modulus one.
For each  $f\in H_E(\xi)$                               
\begin{align*}
(C_{\Upsilon_{1},\Phi_{1}} C_{\Upsilon_{2},\Phi_{2}})f &=C_{\Upsilon_{1},\Phi_{1}}(C_{\Upsilon_{2},\Phi_{2}}f)\\
                                               &=C_{\Upsilon_{1},\Phi_{1}}(\Upsilon_2 \cdot (f \circ \Phi_2 ))\\
                                                &=C_{\Upsilon_{1},\Phi_{1}}(2 (f \circ \Phi_2 ))\\
                                               &=2 (\Upsilon_1 \cdot ( (f \circ \Phi_2 )\circ \Phi_1))\\
                                               &=\big(2 \frac{1}{2} \big) ( f \circ (\Phi_2 \circ \Phi_1))\\
                                               &=C_{\Phi_2 \circ \Phi_1}f. 
\end{align*}                             
Since
$$(\Phi_{2} \circ \Phi_{1})(z)=\Phi_{2}( \Phi_{1}(z))=\Phi_{2}( iz)=-i(iz)=z \,\,\, \text{for all}\,\, z\in \mathbb{C}$$
therefore,
\begin{equation*} 
C_{\Upsilon_{1},\Phi_{1}} C_{\Upsilon_{2},\Phi_{2}} =C_{\Phi_2 \circ \Phi_1}=I_{H_{E}(\xi)}.
\end{equation*}
Hence, $C_{\Upsilon_{1},\Phi_{1}}C_{\Upsilon_{2},\Phi_{2}}$ is unitary and so isometry on $H_{E}(\xi)$.
\end{example}
\begin{theorem} \label{Weigh Com Op isometry upsilon is kernel func}
Let $\Upsilon(z)=K_{c}(z)$ and  $\Phi(z)=\mu z$ for $z\in \mathbb{C}$ where $ \mu,c \in \mathbb{C}$ with $|\mu|\leq 1$ such that $C_{\Upsilon,\Phi}$ is bounded on $H_{E}(\xi).$  Then, $C_{\Upsilon,\Phi}$ is an isometry on $H_{E}(\xi)$ if and only if $|\mu|=1, c=0$ and ${\xi_0}=1.$ Moreover, $C_{\Upsilon,\Phi}$ is unitary on $H_{E}(\xi)$ if and only if $|\mu|=1, c=0$ and ${\xi_0}=1.$
\end{theorem} 
\textbf{Proof.} Let $C_{\Upsilon,\Phi}$ be an isometry. Then, $\mu\neq 0$ and for $g,h\in H_E(\xi)$
\begin{equation*} 
 \langle C_{\Upsilon,\Phi}h, C_{\Upsilon,\Phi}g \rangle=\langle h, g\rangle 
\end{equation*}
that is,
\begin{equation} \label{weigh comp isometry Kernel condition 2 equ}
\langle \Upsilon\cdot (h \circ \Phi), \Upsilon\cdot (g \circ \Phi) \rangle=\langle h, g\rangle.
\end{equation} 
Now we consider the following cases for $h,g \in H_E(\xi)$ :\\
Case I. If $h(z)=1$  for all $z\in \mathbb{C},$ then 
\begin{equation*} 
(\Upsilon\cdot (h \circ \Phi))(z)=\Upsilon(z)h(\mu z)=\Upsilon(z)= \sum_{n=0}^{\infty} \frac{1}{{\xi_n}^2} \overline{c}^n z^n
\end{equation*}
and taking $g=h$ in equation \eqref{weigh comp isometry Kernel condition 2 equ} we get
\begin{equation}  \label{weigh comp isometry Kernel condition 4 equ}
\sum_{n=0}^{\infty} \frac{|c|^{2n}}{{\xi_n}^2} ={\xi_0}^2.
\end{equation}

Case II. If $h(z)=1$ and $g(z)=z$  for all $z\in \mathbb{C},$ then 
\begin{equation*}  
(\Upsilon\cdot (g \circ \Phi))(z)=\Upsilon(z)g(\mu z)=(\mu z)\Upsilon(z)= \sum_{n=0}^{\infty} \frac{1}{{\xi_n}^2} \mu \bar{c\,}^n z^{n+1}=\sum_{n=1}^{\infty} \frac{1}{{\xi_{n-1}}^2} \mu \bar{c\,}^{n-1} z^{n}
\end{equation*}
and
\begin{equation*} 
(\Upsilon\cdot (h \circ \Phi))(z)= \sum_{n=0}^{\infty} \frac{1}{{\xi_n}^2} \overline{c}^n z^n=\frac{1}{{\xi_0}^2}+\sum_{n=1}^{\infty} \frac{1}{{\xi_n}^2} \overline{c}^n z^n.
\end{equation*}
Substituting $h$ and $g$ functions in equation \eqref{weigh comp isometry Kernel condition 2 equ} we obtain
\begin{align*} 
\sum_{n=1}^{\infty} \Big(\frac{\bar{c\,}^{n}}{{\xi_{n}}^2}\Big) \Big( \frac{\bar{\mu}c^{n-1}}{{\xi_{n-1}}^2} \Big) {\xi_{n}}^2  &=0 \\
(\bar{c}\bar{{\mu}}) \sum_{n=1}^{\infty}  \frac{\bar{c\,}^{n-1} c^{n-1}}{{\xi_{n-1}}^2}    &=0 \\
(\bar{c} \bar{\mu}) \sum_{n=0}^{\infty} \frac{|c|^{2n}}{{\xi_n}^2}  &=0.
\end{align*}
Combining equation \eqref{weigh comp isometry Kernel condition 4 equ} with the above equation we get $\bar{c} \bar{\mu} {\xi_0}^2=0$. This implies $c=0$ because $\mu \neq 0$ and ${\xi_0}> 0.$ Then, $\Upsilon(z)=K_{0}(z)=\frac{1}{{\xi_0}^2}$ so 
$$C_{\Upsilon,\Phi}=\frac{1}{{\xi_0}^2} C_{\Phi} \,\,\, \text{and} \,\,C_{\Upsilon,\Phi}^*=\frac{1}{{\xi_0}^2} C_{\Phi}^*.$$
Since $C_{\Upsilon,\Phi}$ is an isometry, $C_{\Upsilon,\Phi}^*C_{\Upsilon,\Phi}=I_{H_E(\xi)}.$ For each $f$ in ${H_E(\xi)}$ and $z\in \mathbb{C}$
\begin{align}
(C_{\Upsilon,\Phi}^*C_{\Upsilon,\Phi}f)(z)&=f(z) \notag \\
\frac{1}{{\xi_0}^2} \frac{1}{{\xi_0}^2} (C_{\Phi}^*C_{\Phi}f)(z) &=f(z) \notag \\
\frac{1}{{\xi_0}^4} f(\mu \bar{\mu} z) &=f(z). \label{weigh comp isometry Kernel condition 6 equ}
\end{align}
In particular for $f(z)=z$ in ${H_E(\xi)}$ and $z=1$ we get
$$\frac{1}{{\xi_0}^4} |\mu|^2=1, \,\,\text{that is,}\,\, |\mu|={{\xi_0}^2}.$$
Further, from equation \eqref{weigh comp isometry Kernel condition 6 equ} for $f(z)=z+1$ in ${H_E(\xi)}$ we get
\begin{align}
\frac{1}{{\xi_0}^4} (|\mu|^2z+1) &=z+1 \notag \\
z+\frac{1}{{\xi_0}^4} &=z+1. \notag 
\end{align}
Thus, $\xi_0=1$ because $\xi_0>0$ and hence, $|\mu|=\xi_{0}^2=1.$\\
Conversely,  let $|\mu|=1, c=0$ and ${\xi_0}=1.$ Then, $\Upsilon \equiv 1$ and this implies $C_{\Upsilon,\Phi}=C_{\Phi}.$ By Theorem \ref {Com Op unitary iff mod(a)=1}, it follows that $C_{\Upsilon,\Phi}$ is unitary and hence an isometry.\\

In Theorem \ref{Weigh Com Op isometry upsilon is kernel func} if we take $\Phi(z)$ to be $\mu z-d$ where $|\mu| \leq 1,$ then for isometry $C_{\Upsilon,\Phi}$ we only succeed in proving the necessary condition that $c=\frac{d}{\mu}.$\\ 
Case I. If $h\equiv 1,$ then $(C_{\Upsilon,\Phi}h)(z)=\frac{1}{{\xi_0}^2}+\sum_{n=1}^{\infty} \frac{1}{{\xi_n}^2} \overline{c}^n z^n.$  Taking $g$ to be the same as $h$ in equation \eqref{weigh comp isometry Kernel condition 2 equ} we obtain 
\begin{equation} \label{general case isometry weigh comp}
\sum_{n=0}^{\infty} \frac{|c|^{2n}}{{\xi_n}^2} ={\xi_0}^2.
\end{equation}
Case II. If $h(z)=1$ and $g(z)=z$  for all $z\in \mathbb{C},$ then 
\begin{align*}  
(\Upsilon\cdot (g \circ \Phi))(z) &=\Upsilon(z)g(\mu z-d)\\
                                 &=\mu z\Upsilon(z)-d\Upsilon(z) \\
                  &= \sum_{n=0}^{\infty} \mu\frac{1}{{\xi_n}^2} \bar{c\,}^n z^{n+1}-\sum_{n=0}^{\infty} d\frac{1}{{\xi_{n}}^2} \bar{c\,}^{n} z^{n} \\
                  &=\left(\frac{-d}{{\xi_{0}}^2}\right)+\sum_{n=1}^{\infty} \left(\mu\frac{1}{{\xi_{n-1}}^2} \bar{c\,}^{n-1} - d\frac{1}{{\xi_{n}}^2} \bar{c\,}^{n} \right) z^{n}. \\
\end{align*} 
For these $g$ and $h,$ equation \eqref{weigh comp isometry Kernel condition 2 equ} gives
\begin{align*}
 \left(\frac{1}{{\xi_{0}}^2}\right) \left(\frac{-\bar{d}}{{\xi_{0}}^2}\right) {\xi_{0}}^2 +\sum_{n=1}^{\infty} \left(\frac{1}{{\xi_n}^2} \bar{c\,}^n\right) \left( \bar{\mu} \frac{1}{{\xi_{n-1}}^2} {c}^{n-1} - \bar{d}\frac{1}{{\xi_{n}}^2} {c}^{n} \right) {\xi_n}^2 &=0 \\
  \left(\frac{-\bar{d}}{{\xi_{0}}^2}\right)+\sum_{n=1}^{\infty} \left( \bar{\mu} \bar{c} \frac{|c|^{2(n-1)}} {{\xi_{n-1}}^2} -\bar{d} \frac{|c|^{2n}}{{\xi_{n}}^2}  \right) &=0 \\
 \bar{\mu}\bar{c} \sum_{n=0}^{\infty} \frac{|c|^{2n}} {{\xi_{n}}^2} -\bar{d}\left( \frac{1}{{\xi_{0}}^2}+\sum_{n=1}^{\infty} \frac{|c|^{2n}}{{\xi_{n}}^2} \right)=0.
\end{align*} 
Combining equation \eqref{general case isometry weigh comp} with the above equation we get 
\begin{equation}
(\bar{\mu}\bar{c}-\bar{d}){\xi_{0}}^2=0.
\end{equation}
It clearly follows that $\mu \neq 0$ because if $\mu= 0,$ then $d=0$ and $C_{\Upsilon,\Phi}$ will not be isometry. Hence, $c=\frac{d}{\mu}.$

\section{$m$-isometry and invertible composition operator} 
Let $H$ be a Hilbert space, m be a positive integer and $S$ be a bounded linear operator on $H$ then $S$ is called $m$-isometry (see \cite{m-isometry Agler}) if it satisfies
$$\sum_{k=0}^{m} (-1)^{m-k}\,\, {m \choose k} {S^*}^k S^k=0.$$
In the following result we characterize $m$-isometry composition operator in terms of isometry on $H_E(\xi)$.
\begin{theorem}
Let $\mu\neq 0$ with $|\mu|\leq 1$ and $\Phi(z)=\mu z$ for $z\in \mathbb{C}.$ Then $C_{\Phi}$ is m-isometry on $H_E(\xi)$ if and only if $C_{\Phi}$ is isometry.
\end{theorem}
\textbf{Proof.} Let $\Psi(z)=\overline{\mu}z$ for all $z\in \mathbb{C}.$ Then $C_{\Phi}^*=C_{\Psi},$ therefore
\begin{align*} 
\sum_{k=0}^{m} ((-1)^{m-k}\,\, {m \choose k} {C_{\Phi}^*}^k C_{\Phi}^kf)(z)&=\sum_{k=0}^{m} ((-1)^{m-k}\,\, {m \choose k} C_{\Psi}^k C_{\Phi}^kf)(z)\\
&=\sum_{k=0}^{m} (-1)^{m-k}\,\, {m \choose k} \,f(|\mu|^{2k} z). 
\end{align*}
Let $C_{\Phi}$ be $m$-isometry. Then $\sum_{k=0}^{m} (-1)^{m-k}\,\, {m \choose k} {C_{\Phi}^*}^k C_{\Phi}^kf=0$ for all $f$ in $H_E(\xi).$ In particular, $f$ in $H_E(\xi),$ defined by $f(z)=z$ for all $z\in \mathbb{C}$ gives 
$$\sum_{k=0}^{m} (-1)^{m-k}\,\, {m \choose k} \,|\mu|^{2k}z =0.$$
Further, for $z=1$ this implies
$\sum_{k=0}^{m} (-1)^{m-k}\,\, {m \choose k} \,|\mu|^{2k} =0,$ that is, $(1-|\mu|^2)^m=0.$ Hence, $|\mu|=1$ and by Theorem \ref{Com Op isometry iff mod(a)=1}, $C_{\Phi}$ is isometry.

Conversely, let $C_{\Phi}$ be isometry. Then, $C_{\Phi}^*C_{\Phi}=I_{H_E(\xi)}$ which further implies $C_{\Phi}$ is $m$-isometry.\\
\begin{proposition} \label{product of two comp op calculated}
Let $C_{\Upsilon_{1},\Phi_{1}}$ and $C_{\Upsilon_{2},\Phi_{2}}$ be two bounded linear operators on $H_E(\xi)$ then 
$$C_{\Upsilon_{1},\Phi_{1}} C_{\Upsilon_{2},\Phi_{2}} =C_{\Upsilon_{1} (\Upsilon_{2} \circ \Phi_{1}),\Phi_{2} \circ \Phi_{1}}.$$
\end{proposition}
\textbf{Proof.} For each  $f\in H_E(\xi)$                               
\begin{align*}
(C_{\Upsilon_{1},\Phi_{1}} C_{\Upsilon_{2},\Phi_{2}})f &=C_{\Upsilon_{1},\Phi_{1}}(C_{\Upsilon_{2},\Phi_{2}}f)\\
                                               &=C_{\Upsilon_{1},\Phi_{1}}(\Upsilon_2\cdot (f \circ \Phi_2 ))\\
                                               &=\Upsilon_1 \cdot ((\Upsilon_2 \cdot (f \circ \Phi_2 ))\circ \Phi_1) 
\end{align*}
and $z\in \mathbb{C}$
\begin{align*}
((C_{\Upsilon_{1},\Phi_{1}} C_{\Upsilon_{2},\Phi_{2}})f)(z)&=\Upsilon_1(z) ((\Upsilon_2 \cdot (f \circ \Phi_2 ))\circ \Phi_1)(z)\\
                                                   &=\Upsilon_1(z) (\Upsilon_2 \cdot (f \circ \Phi_2 )) ( \Phi_1(z))\\
                                                    &=\Upsilon_1(z) \Upsilon_2( \Phi_1(z)) (f \circ \Phi_2 ) ( \Phi_1(z))\\
                                                     &=\Upsilon_1(z) \Upsilon_2( \Phi_1(z)) f( \Phi_2 ( \phi_1(z)))\\ 
                                                    &=\Upsilon_1(z) (\Upsilon_2 \circ \Phi_1)(z) (f \circ (\Phi_2 \circ \Phi_1))(z)\\                                                    
                                                    &=( C_{\Upsilon_{1} (\Upsilon_{2} \circ \Phi_{1}),\Phi_{2} \circ \Phi_{1}}f)(z).
\end{align*}
Thus, $C_{\Upsilon_{1},\Phi_{1}} C_{\Upsilon_{2},\Phi_{2}} =C_{\Upsilon_{1} (\Upsilon_{2} \circ \Phi_{1}),\Phi_{2} \circ \Phi_{1}}.$\\ 

For positive integer $k,$ let $\Phi_{[k]}$ denote the  $k$ times composition of $\Phi$ with itself and in particular $\Phi_{[1]}=\Phi.$ From Proposition \ref{product of two comp op calculated}, it follows that
$C_{\Upsilon,\Phi}^2=C_{\Upsilon,\Phi} C_{\Upsilon,\Phi} =C_{\Upsilon (\Upsilon \circ \Phi_{[1]}),\Phi_{[2]}}.$ Further, by applying the Principle of Mathematical Induction the following result follows:
\begin{corollary} \label{weighted comp op power k}
Let $k>1$ be a positive integer and $C_{\Upsilon,\Phi}$ be bounded on $H_E(\xi).$ Then, 
\begin{equation*} 
C_{\Upsilon,\Phi}^k =C_{\Upsilon (\Upsilon \circ \Phi_{[1]})(\Upsilon \circ \Phi_{[2]})...(\Upsilon \circ \Phi_{[k-1]}), \Phi_{[k]}}.
\end{equation*}
\end{corollary}
\begin{theorem}
Let $\Upsilon(z)=\nu z$ and $\Phi(z)=\mu z$  for $z\in \mathbb{C}$ where $\nu,\mu \in \mathbb{C}$ with $|\mu|\leq 1$  such that $C_{\Upsilon,\Phi}$ is bounded on $H_{E}(\xi).$ Then, $C_{\Upsilon,\Phi}$ is $m$-isometry if and only if 
\begin{equation*} 
\sum_{k=1}^{m} (-1)^{m-k}\,\, {m \choose k} {\xi_{n+k}}^{2} |\nu|^{2k}  |\mu|^{k(2n+k-1)}=(-1)^{m+1} {\xi_{n}^2} \,\,\, \text{for all} \,\,\,n \geq 0.
\end{equation*}
\end{theorem}
\textbf{Proof.} Let $g(z)=\sum_{n=0}^{\infty} b_{n} z_{n} \in H_{E}(\xi).$ Then, using Corollary \ref{weighted comp op power k} for $k>1$ we get
\begin{align*} 
(C_{\Upsilon,\Phi}^k g)(z) &=(C_{\Upsilon (\Upsilon \circ \Phi_{[1]})(\Upsilon \circ \Phi_{[2]})...(\Upsilon \circ \Phi_{[k-1]}), \Phi_{[k]}}g)(z) \\
                           &=\Upsilon(z) (\Upsilon \circ \Phi_{[1]})(z)(\Upsilon \circ \Phi_{[2]})(z)...(\Upsilon \circ \Phi_{[k-1]})(z) g(\Phi_{[k]}(z)) \\
                           &=\nu^{k} \mu^{1+2+3+...(k-1)} z^{k} \sum_{n=0}^{\infty} b_{n} (\mu^{k} z)^n \\
                            &= \sum_{n=0}^{\infty} b_{n} \nu^{k} \mu^{\frac{k(2n+k-1)}{2}}  z^{n+k} \\
                           &= \sum_{n=k}^{\infty} b_{n-k} \nu^{k} \mu^{\frac{k(2n-k-1)}{2}} z^{n}.                      
\end{align*}
The above calculation also holds for $k=1$ since $C_{\Upsilon,\Phi}^1=C_{\Upsilon,\Phi},$ therefore,
$$(C_{\Upsilon,\Phi}^1g)(z)=(C_{\Upsilon,\Phi}g)(z)=\sum_{n=0}^{\infty} b_{n} \nu \mu^{n}  z^{n+1}= \sum_{n=1}^{\infty} b_{n-1} \nu \mu^{n-1} z^{n}.$$
This implies for $K_{w}(z)=\sum_{n=0}^{\infty} \frac{1}{\xi_{n}^2} \bar{w}^n z^{n}$ and $e_{n}(z)=\frac{1}{\xi_{n}}z^n$ 
\begin{equation} \label{weigh comp value at rep kerne}
(C_{\Upsilon,\Phi}^k K_{w})(z)=\sum_{n=k}^{\infty} \frac{1}{\xi_{n-k}^2} \bar{w}^{n-k} \nu^{k} \mu^{\frac{k(2n-k-1)}{2}} z^{n}                      
\end{equation}
and
\begin{align}
(C_{\Upsilon,\Phi}^k e_{n})(z) &=\nu^{k} \mu^{1+2+3+...(k-1)} z^{k} \frac{1}{\xi_{n}} (\mu^{k} z)^n \notag \\
                               &=\frac{1}{\xi_{n}} \nu^{k} \mu^{\frac{k(2n+k-1)}{2}}   z^{n+k}. \label{wei comp power k evaluated at basis}
\end{align}
Let $((C_{\Upsilon,\Phi}^k)^{*} g)(z)=\sum_{n=0}^{\infty} r_{n}z^{n}$ then, 
\begin{align*}
\langle (C_{\Upsilon,\Phi}^k)^{*} g, e_{n} \rangle &= r_{n} \xi_{n}.
\end{align*}
In fact using \eqref{wei comp power k evaluated at basis} we get
\begin{align*}
r_{n} &=\frac{1}{\xi_n} \langle (C_{\Upsilon,\Phi}^k)^{*} g, e_{n} \rangle \\ 
       &=\frac{1}{\xi_n} \langle g, C_{\Upsilon,\Phi}^k e_{n} \rangle \\
        &=\frac{1}{\xi_n}  b_{n+k}\left(\frac{1}{\xi_{n}} \bar{\nu}^{k} \bar{\mu}^{\frac{k(2n+k-1)}{2}}  \right) {\xi_{n+k}}^{2}\\ 
       &= \frac{{\xi_{n+k}}^{2}}{\xi_n^{2}}  b_{n+k} \bar{\nu}^{k} \bar{\mu}^{\frac{k(2n+k-1)}{2}} . 
\end{align*}
Hence 
\begin{equation} \label{weigh comp adjoint value in general}
((C_{\Upsilon,\Phi}^k)^{*} g)(z) =\sum_{n=0}^{\infty} \left(\frac{{\xi_{n+k}}^{2}}{\xi_n^{2}}  b_{n+k} \bar{\nu}^{k} \bar{\mu}^{\frac{k(2n+k-1)}{2}}  \right) z^{n}.
\end{equation}
Applying equations \eqref{weigh comp value at rep kerne} and \eqref{weigh comp adjoint value in general} we get
\begin{align}
((C_{\Upsilon,\Phi}^k)^{*}C_{\Upsilon,\Phi}^{k} K_{w})(z) &=\sum_{n=0}^{\infty} \left(\frac{{\xi_{n+k}}^{2}}{\xi_n^{2}}  \Big(\frac{1}{\xi_{n}^2} \bar{w}^{n} \nu^{k} \mu^{\frac{k(2n+k-1)}{2}}   \Big) \bar{\nu}^{k} \bar{\mu}^{\frac{k(2n+k-1)}{2}}  \right) z^{n} \notag\\
                                 &= \sum_{n=0}^{\infty} \frac{{\xi_{n+k}}^{2}}{\xi_n^{2} \, \xi_n^{2}} |\nu|^{2k}  |\mu|^{k(2n+k-1)} \bar{w}^{n} z^{n}. \label{product power k calculated for rep kernel}
\end{align}
It is well known that for a bounded operator $S,$ ${S^*}^k=(S^k)^{*}.$ This along with \eqref{product power k calculated for rep kernel} gives  
\begin{align*}
\sum_{k=1}^{m} (-1)^{m-k}\,\, {m \choose k} ({C_{\Upsilon,\Phi}^*}^k {C_{\Upsilon,\Phi}}^k K_{w})(z) &=\sum_{k=1}^{m} (-1)^{m-k}\,\, {m \choose k} \Big(\sum_{n=0}^{\infty} \frac{{\xi_{n+k}}^{2}}{\xi_n^{2} \, \xi_n^{2}} |\nu|^{2k}  |\mu|^{k(2n+k-1)} \bar{w}^{n} z^{n} \Big) \\
                              &=\sum_{n=0}^{\infty} \Big(\sum_{k=1}^{m} (-1)^{m-k}\,\, {m \choose k} \frac{{\xi_{n+k}}^{2}}{\xi_n^{2} \, \xi_n^{2}} |\nu|^{2k}  |\mu|^{k(2n+k-1)} \Big) \bar{w}^{n} z^{n}. 
\end{align*}
Since $C_{\Upsilon,\Phi}$ is bounded on $H_{E}(\xi)$ and set $Span\{K_{w} : w\in \mathbb{C} \}$ is dense in $H_{E}(\xi),$  therefore, $C_{\Upsilon,\Phi}$ is $m$-isometry if and only if
\begin{equation*}
\sum_{k=0}^{m} (-1)^{m-k}\,\, {m \choose k} ({C_{\Upsilon,\Phi}^*}^k {C_{\Upsilon,\Phi}}^k K_{w})(z)=0 \,\,\, \text{for all}\,\, z, w \in \mathbb{C}
\end{equation*}
equivalently,
\begin{align*}
\sum_{k=1}^{m} (-1)^{m-k}\,\, {m \choose k} ({C_{\Upsilon,\Phi}^*}^k {C_{\Upsilon,\Phi}}^k K_{w})(z) &=(-1)^{m+1} (I_{H_{E}(\xi)}K_{w})(z) \\ 
\sum_{n=0}^{\infty} \Big(\sum_{k=1}^{m} (-1)^{m-k}\,\, {m \choose k} \frac{{\xi_{n+k}}^{2}}{\xi_n^{2} \, \xi_n^{2}} |\nu|^{2k}  |\mu|^{k(2n+k-1)} \Big) \bar{w}^{n} z^{n} &=\sum_{n=0}^{\infty} (-1)^{m+1} \frac{1}{\xi_{n}^2} \bar{w}^n z^{n}.
\end{align*}
Comparing coefficients on both sides we get $C_{\Upsilon,\Phi}$ is $m$-isometry if and only if
\begin{equation*} 
\sum_{k=1}^{m} (-1)^{m-k}\,\, {m \choose k} {\xi_{n+k}}^{2} |\nu|^{2k}  |\mu|^{k(2n+k-1)}=(-1)^{m+1} {\xi_{n}^2} \,\,\, \text{for all} \,\,n \geq 0.
\end{equation*}\\

For $\mu (\neq 0) \in \mathbb{C}$ and for sequence $\mu\xi :=\big(\frac{1}{|\mu|^n}\xi_n\big)_{n\in \mathbb{Z_+}}$ we have
$$lim_{n\rightarrow \infty} \left(\frac{1}{|\mu|^n} \xi_n \right)^{1/n}=\frac{1}{|\mu|} lim_{n\rightarrow \infty} \xi_n ^{1/n} =\infty$$
and Hilbert space of entire functions $H_E(\mu\xi)$ is  
$$H_E(\mu\xi)=\Big\{ f(z)=\sum_{n=0}^{\infty} b_n z^n : ||f||^2=\sum_{n=0}^{\infty} |b_n|^2 \frac{1}{|\mu|^{2n}}{\xi_n}^2 <\infty \Big\}.$$ 

It is not difficult to see that for $\mu\neq 0$ with $|\mu|\leq 1,$ function $f(z)=\sum_{n=0}^{\infty} b_n z^n$ in $H_E(\mu\xi )$ satisfies
$$\sum_{n=0}^{\infty} |b_n|^2 {\xi_n}^2  \leq \sum_{n=0}^{\infty} |b_n|^2 \frac{1}{|\mu|^{2n}}{\xi_n}^2 <\infty. $$
Thus, $H_E(\mu\xi)\subseteq H_E(\xi)$ and for $|\mu|=1,$ $H_E(\mu\xi)= H_E(\xi).$ The following result deals with the invertibility.
\begin{theorem}
Let $\mu\neq 0$ with $|\mu|\leq 1$ and let functions $\Phi$ and $\Psi$ be defined by $\Phi(z)=\mu z$ and $\Psi(z)=\frac{1}{\mu} z$ for $z\in \mathbb{C}.$ Then  $C_{\Phi}:H_E(\xi) \longrightarrow H_E(\mu\xi )$ is invertible with inverse $C_{\Psi}:H_E(\mu\xi) \longrightarrow H_E(\xi). $ 
\end{theorem} 
\textbf{Proof.} Clearly, if $f\in H_E(\mu\xi )$ and $g\in H_E(\xi)$ then $C_{\Psi}f \in H_E(\xi)$ and $C_{\Phi}g \in H_E(\mu\xi ).$ Now, for $f\in H_E(\mu\xi),$ $g\in H_E(\xi)$ and $z\in \mathbb{C}$ 
$$(C_{\Phi} C_{\Psi}f)(z)= f(\Psi(\Phi(z)))=f(z)$$
$$(C_{\Psi} C_{\Phi}g)(z)= g(\Phi(\Psi(z)))=g(z).$$

Thus, $C_{\Phi} C_{\Psi}=I_{H_E(\mu\xi )}$ and $C_{\Psi} C_{\Phi}=I_{H_E(\xi)}.$  Hence, $C_{\Phi}$ is invertible and $C_{\Phi}^{-1}=C_{\Psi}.$ 
\begin{example}
Let $\Phi(z)=(\frac{1}{2}+\frac{1}{2}i) z$ and $\Psi(z)=\frac{1}{1-i} z$ for $z\in \mathbb{C}.$  Then, operator $C_{\Psi}:H_E((\frac{1}{2}+\frac{1}{2}i)\xi) \longrightarrow H_E(\xi)$ is inverse of
$C_{\Phi}:H_E(\xi) \longrightarrow H_E((\frac{1}{2}+\frac{1}{2}i)\xi )$ 
\end{example}
\section{Generalized weighted composition operator $\mathcal{D}_{\Phi,\Upsilon}^p$}
Kumar and  Abbas \cite{Pawan and Zaheer} obtained the necessary and sufficient condition for operator $\mathcal{D}_{\Phi,\Upsilon}^{1}$ to be bounded on weighted Hardy spaces of holomorphic functions on $\mathbb{D},$ the open unit disk.
On the similar lines we generalize this for operator $\mathcal{D}_{\Phi,\Upsilon}^p$ on $H_E(\xi)$ and prove the following  result.
\begin{theorem} \label{generalized weighted com op bounded thm}
Let $\Phi,\Upsilon$ be entire functions on $ \mathbb{C}$ such that $ \{\Phi^n\}_{n\in \mathbb{Z_+}}$ is orthogonal. Then the operator $\mathcal{D}_{\Phi,\Upsilon}^p$ is bounded on $H_E(\xi)$ if and only if 
$$||\Upsilon \Phi^{n-p}||\leq \frac{M \xi_n}{n(n-1)(n-2)...(n-(p-1))} \,\,\, \text{for all}\,\, n\geq p.$$
\end{theorem}
\textbf{Proof.} 
Now for $n\geq p$
\begin{align}
 ||\mathcal{D}_{\Phi,\Upsilon}^p z^n||&=n||\mathcal{D}_{\Phi,\Upsilon}^{p-1}(z^{n-1})|| \notag\\ 
                               &=n(n-1)||\mathcal{D}_{\Phi,\Upsilon}^{p-2}(z^{n-2})|| \notag\\ 
                             &=n(n- 1)(n-2)...(n-(p-1))||\Upsilon \Phi^{n-p}||. \label{equation bdd norm thm1 eq 1}
\end{align}
Now let $\mathcal{D}_{\Phi,\Upsilon}^p$ be a bounded linear operator then, there exists a number $M$ such that 
\begin{equation} \label{bounded thm1 equ have M}
||\mathcal{D}_{\Phi,\Upsilon}^p f|| \leq M ||f|| \,\,\, \text{for all}\, f\in Dom(\mathcal{D}_{\Phi,\Upsilon}^p) .
\end{equation}
 From equations (\ref{equation bdd norm thm1 eq 1}) and (\ref{bounded thm1 equ have M}) it follows that for $n\geq p$
 \begin{equation*}
 n(n- 1)(n-2)...(n-(p-1))||\Upsilon  \Phi^{n-p}|| \leq M ||z^n||=M \xi_n 
 \end{equation*}
that is, $||\Upsilon  \Phi^{n-p}||\leq \frac{M \xi_n}{n(n-1)(n-2)...(n-(p-1))}.  $\\
 Conversely, let for $n\geq p$  
 \begin{equation} \label{converse let satisfy condition equation}
 ||\Upsilon \Phi^{n-p}||\leq \frac{M \xi_n }{n(n-1)(n-2)...(n-(p-1))}. 
 \end{equation}
 For $f(z)=\sum_{n=0}^{\infty} b_n z^n \in Dom(\mathcal{D}_{\Phi,\Upsilon}^p),$
 $$f^p(z)=\sum_{n=p}^{\infty} n(n-1)(n-2)...(n-(p-1)) b_n z^{n-p}.$$ Now using the fact that $ \{\Phi^n\}_{n\geq 0}$ is orthogonal we have 
 \begin{align*}
||\mathcal{D}_{\Phi,\Upsilon}^pf||^2&=||\Upsilon (f^p \circ \Phi)||^2\\
                   &=||\Upsilon  \sum_{n=p}^{\infty} n(n-1)(n-2)...(n-(p-1)) b_n \Phi^{n-p}||^2\\
                   &=|| \sum_{n=p}^{\infty} n(n-1)(n-2)...(n-(p-1)) b_n\, \Upsilon \Phi^{n-p}||^2\\
                   &=\sum_{n=p}^{\infty} ||n(n-1)(n-2)...(n-(p-1)) b_n\, \Upsilon \Phi^{n-p}||^2\\
                   &=\sum_{n=p}^{\infty} n^2(n-1)^2(n-2)^2...(n-(p-1))^2 |b_n|^2\,||\Upsilon \Phi^{n-p}||^2
\end{align*}
Applying equation \eqref{converse let satisfy condition equation} we get
\begin{align*}
||\mathcal{D}_{\Phi,\Upsilon}^pf||^2 &\leq\sum_{n=p}^{\infty} n^2(n-1)^2(n-2)^2...(n-(p-1))^2 |b_n|^2\, \frac{M^2 \,\xi_n^2}{n^2(n-1)^2(n-2)^2...(n-(p-1))^2}\\
                    &=M^2 \sum_{n=p}^{\infty} |b_n|^2  \xi_n^2 \\  
                    & \leq M^2 \sum_{n=0}^{\infty} |b_n|^2  \xi_n^2\\              
                    &=M^2 ||f||^2.
 \end{align*}
 Hence $||\mathcal{D}_{\Phi,\Upsilon}^pf|| \leq M||f||$ for all $f \in Dom(\mathcal{D}_{\Phi,\Upsilon}^p).$\\


If $\Upsilon \equiv 1$ then the Theorem \ref{generalized weighted com op bounded thm} reduces to the following result:
\begin{corollary}
Let $\Phi$ be entire function on $ \mathbb{C}$ such that $ \{\Phi^n\}_{n\in \mathbb{Z_+}}$ is orthogonal. Then the operator $C_{\Phi}  \mathcal{D}^p$ is bounded on $H_E(\xi)$ if and only if 
$$|| \Phi^{n-p}||\leq \frac{M \xi_n}{n(n-1)(n-2)...(n-(p-1))} \,\, \text{for all}\,\, n\geq p.$$\\
\end{corollary}

If $\Upsilon \equiv {\Phi}'$ and  $p=1$ then $\mathcal{D}_{\Phi,\Upsilon}^pf=\Phi' (f'\circ \Phi)=(f\circ \Phi)'=   \mathcal{D}C_{\Phi}f.$ Combining this with Theorem \ref{generalized weighted com op bounded thm} we have the following result:

\begin{corollary}
Let $\Phi$ be entire function on $ \mathbb{C}$ such that $ \{\Phi^n\}_{n\in \mathbb{Z_+}}$ is orthogonal. Then, the operator $\mathcal{D}C_{\Phi}$ is bounded on $H_E(\xi)$ if and only if 
$$|| {\Phi}'\Phi^{n-1}||\leq \frac{M \xi_n}{n}.$$
\end{corollary}

\textbf{Anuradha Gupta}\\
 Department of Mathematics, Delhi College of Arts and Commerce,\\
  University of Delhi, New Delhi-110023, India.\\
  \vspace{0.2cm}
 email: dishna2@yahoo.in\\
  \textbf{Geeta Yadav}\\
  Department of Mathematics, University of Delhi, \\
  New Delhi-110007, India.\\
  email: ageetayadav@gmail.com
\end{document}